\documentclass[11 pt]{amsart}
\usepackage{amsmath}
\usepackage{amssymb}
\usepackage{amsthm}
\usepackage{epsfig}
\usepackage[margin=1.25in]{geometry}
\usepackage{amsfonts}
\usepackage[font=small,format=plain,labelfont=bf,up,textfont=it,up]{caption}
\usepackage{amscd}
\usepackage[all,cmtip]{xy}

\newtheorem{theorem}{Theorem}[section]

\newtheorem{lemma}[theorem]{Lemma}
\newtheorem{corollary}[theorem]{Corollary}

\newtheorem*{namedtheorem}{\theoremname}
\newcommand{\theoremname}{testing}
\newenvironment{named}[1]{\renewcommand{\theoremname}{#1}\begin{namedtheorem}}{\end{namedtheorem}}

\theoremstyle{definition}

\theoremstyle{remark}
\newtheorem*{remark}{Remark}

\newtheorem*{acknowledgements}{Acknowledgments}

\DeclareMathOperator{\Homeo}{Homeo}
\DeclareMathOperator{\Ker}{ker}

\DeclareMathOperator{\Mod}{Mod}
\newcommand\Torelli{\text{${\mathcal I}$}}
\DeclareMathOperator{\Sp}{Sp}
\DeclareMathOperator{\PSp}{PSp}

\DeclareMathOperator{\interior}{Int}

\newcommand\R{\text{$\mathbb{R}$}}
\newcommand\C{\text{$\mathbb{C}$}}
\newcommand\Z{\text{$\mathbb{Z}$}}
\newcommand\Q{\text{$\mathbb{Q}$}}

\DeclareMathOperator{\HH}{H}

\DeclareMathOperator{\Aut}{Aut}

\DeclareMathOperator{\Id}{Id}

\begin{document}
\thanks{The second author has been partially supported by the NSF grant DMS-0706878, NSF CAREER award 0952106 and the Alfred P. Sloan Foundation}
\begin{abstract}
We show that, for any (symmetric) finite generating set of the Torelli group of a closed surface, the probability that a random word is not pseudo-Anosov decays exponentially
in the length of the word. 
\end{abstract}
\title{On Genericity of Pseudo-Anosovs in the Torelli Group}

\author{Justin Malestein, Juan Souto}

\maketitle
\section{Introduction}
A well-known and fundamental result in the theory of mapping class groups is the Nielsen-Thurston classification, which states that any mapping class must be periodic,
reducible, or pseudo-Anosov.  The mapping class group $\Mod(\Sigma)$ of a surface $\Sigma$ is the group of orientation-preserving homeomorphisms modulo those isotopic to the identity.
Periodic mapping classes are those of finite order and reducible mapping classes are those that permute a collection of disjoint simple closed curves.  Pseudo-Anosov mapping classes are those classes which have a representative 
preserving a pair of transverse projective measured foliations.  A natural question is: which kind of mapping class (if any) is ``typical'' or ``generic?''  

In \cite{Riv2, Riv}, Rivin showed that the set of pseudo-Anosov elements is generic in the mapping class group $\Mod(\Sigma)$. More precisely, he proved that the probability of a random walk (in a symmetric generating set) not being pseudo-Anosov decays exponentially in the length of the word.  (See also \cite{Kowalski} for a proof of this fact by Kowalski using the method of large sieves.) In this paper, we adapt Rivin's arguments to prove that the same result holds true for the Torelli subgroup $\Torelli(\Sigma)$ of the mapping class group $\Mod(\Sigma)$ of a closed surface $\Sigma$ of at least genus $3$.

Before stating the theorem, we set up some notation for random walks.  For a symmetric finite generating set $S$ of a group $G$, we define $W_n(S)$ to be the set of all words, not necessarily
reduced, in $S$ of length $n$.  Let $h(w)\in G$ be the group element corresponding to the word $w \in W_n(S)$.  The following theorem implies that, in particular, pseudo-Anosov mapping
classes are generic in the Torelli group in the sense of random walks, and in fact the complimentary set is exponentially thin.

\begin{theorem} \label{theorem:maintheorem}  
Let $\Sigma$ be a closed surface of genus $g \geq 3$, let $S$ a finite symmetric generating set of $\Torelli(\Sigma)$, and denote by $P_n$ the probability that the element $h(w)\in\Torelli(\Sigma)$ represented by a word $w \in W_n(S)$ chosen uniformly at random is pseudo-Anosov.  Then, $1- P_n = O(c^n)$ for some $c < 1$.
\end{theorem}

Recall that the Torelli group $\Torelli(\Sigma)$ is the subgroup of $\Mod(\Sigma)$ consisting of those elements which act trivially on the integer homology of $\Sigma$; it is due to Johnson \cite{Johnson} that as long as $\Sigma$ has at least genus $3$, the Torelli group is finitely generated. On the other hand, the Torelli subgroup of the mapping class group is not finitely generated in genus $2$ \cite{MM} although some of the results we prove still hold in that case.

\begin{remark}
Before going any further, it should be remarked that it was already known the set of pseudo-Anosov elements is generic in the Torelli group in the sense that 
$\displaystyle \lim_{n \to \infty} P_n = 1$. In fact, Maher \cite{Maher} proved that this is the case in every finitely generated subgroup of the mapping class group which contains some pseudo-Anosov element. Maher's arguments are very different from both Rivin's \cite{Riv2,Riv} and from those in this note and we do not know if they can be used to give a different proof of Theorem \ref{theorem:maintheorem}.
\end{remark}

We outline briefly the proof of Theorem \ref{theorem:maintheorem}. We associate to every 2-sheeted cover $p:\hat\Sigma\to\Sigma$ a homomorphism $\ell_p: \Gamma_p\to \PSp(2g-2, \Z)$ where $\Gamma_p$ is the finite index subgroup of $\Mod(\Sigma)$ consisting of mapping classes which lift to the cover.  The subgroup $\Gamma_p$ contains the Torelli group $\Torelli(\Sigma)$, and using Margulis's Normal Subgroup Theorem we conclude that $\ell_p(\Torelli(\Sigma))$ has finite index in $\PSp(2g-2, \Z)$. Then, we prove a Casson-like criterion: if for every 2-sheeted cover as above, the image of $\varphi\in\Torelli(\Sigma)$ under $\ell_p$ is represented by an element in $\Sp(2g-2,\Z)$ with irreducible characteristic polynomial, then $\varphi$ is pseudo-Anosov. The final ingredient of the proof of Theorem \ref{theorem:maintheorem} is a result due to Rivin \cite{Rivin} asserting that the probability of having reducible characteristic polynomial decays exponentially under a random walk in a Zariski-dense subgroup of $\Sp(2g-2, \Z)$. 

Before moving on we would like to remark that with a small modification of our methods, the same theorem holds for a larger class of groups which are contained in the intersection of
all $\Gamma_p$.  In particular, we can show the following.

\begin{theorem} \label{theorem:maintheorem2}
Let $\Sigma$ be a surface of genus $g > 2$.  Let $G$ be a finitely generated subgroup of $\cap_p \Gamma_p$ with symmetric generating set $S$, and let $G_p < \Sp(2g-2, \Z)$ be the full preimage of $\ell_p(G) < \PSp(2g-2, \Z)$.  If $G_p$ is Zariski dense for all index $2$ covers $p$, 
then the probability of a random walk in $S$ being pseudo-Anosov is $1-O(c^n)$ for some $c < 1$.
\end{theorem}

Now, let ${\mathcal K}(\Sigma)$ be the subgroup of the Torelli group $\Torelli(\Sigma)$ generated by Dehn twists about separating curves. We obtain the following as a corollary to Theorem
\ref{theorem:maintheorem2} and Corollary \ref{corollary:Kgfi}.

\begin{corollary} \label{cor:maincor}
Let $\Sigma$ be closed surface of genus $g > 2$, and let $S$ be a finite (symmetric) generating set of a group $\Lambda$ with ${\mathcal K}(\Sigma) < \Lambda < \cap_p \Gamma_p$, and let $P_n$ be the probability that a 
random word of length $n$ in $S$ is pseudo-Anosov.  Then $1-P_n = O(c^n)$ for some $c < 1$.
\end{corollary}
This, in particular, implies that pseudo-Anosov mapping classes are generic in any finite index subgroup of $\Torelli(\Sigma)$.

\begin{remark}
Lubotzky and Meiri \cite{Lubotzky} have also proven Theorem \ref{theorem:maintheorem} via the same general strategy of proof.  These authors and ourselves have worked completely independently of each other, and we all learned about the other's results once our respective works were almost completed.
\end{remark}

\begin{acknowledgements}
The second author is very thankful for the hospitality of Temple University and the University of Pennsylvania.  The first author would like to thank Igor Rivin for comments on previous drafts and discussions of his results relevant to this paper  and Ben McReynolds for
some helpful discussions on Zariski density.
\end{acknowledgements}

\section{A Few Preliminaries} \label{section:prelim}
We recall a few well-known facts about the mapping class group $\Mod(\Sigma)$ and Torelli group $\Torelli(\Sigma)$ of a closed oriented surface $\Sigma$.  Most of the following definitions
and results may be found in \cite{FarbMarg}.

\subsection{Homology}
In the following, the first homology group of the surface $\Sigma$ with coefficients in $R=\Z,\Z/2\Z$ will routinely be used. Hoping that no confusion will occur, we denote elements in $\pi_1(\Sigma)$ and their class in $\HH_1(\Sigma,R)$ by the same symbol. We will also denote the algebraic intersection pairing for $R=\Z$ and for $R=\Z/2\Z$ in the same way:
$$\langle\cdot,\cdot\rangle:\HH_1(\Sigma,R)\times\HH_1(\Sigma,R)\to R$$
It is well-known that the intersection form $\langle\cdot,\cdot\rangle$ is a non-degenerate alternating form. In fact, if $\Sigma$ has genus $g$ we have 
\begin{equation}\label{eq:preli1}
(\HH_1(\Sigma,R),\langle\cdot,\cdot\rangle)\cong(R^{2g},\langle\cdot,\cdot\rangle)
\end{equation}
where the second bracket is the standard symplectic form.  (See \cite[Section 6.1]{FarbMarg} for a discussion of the symplectic structure.) In particular, any choice of an identification as in \eqref{eq:preli1} yields an identification
\begin{equation}\label{eq:preli2}
\Aut(\HH_1(\Sigma,R),\langle\cdot,\cdot\rangle)\cong\Sp(2g,R)
\end{equation}
Continuing with the same notation, recall that the homomorphisms
$$\Mod(\Sigma)\to\Sp(2g,\Z),\ \ \Mod(\Sigma)\to\Sp(2g,\Z/2\Z)$$
are surjective \cite[Theorem 6.4]{FarbMarg}. Recall that, by definition the Torelli group $\Torelli(\Sigma)$ is the kernel of the first of these two homomorphisms.

Before moving on observe that everything we just said holds also for compact surfaces with a single boundary component.

\subsection{Covers and lifts}
Below we will consider covers $p:\hat\Sigma\to\Sigma$ of degree $2$ with $\hat\Sigma$ connected. It is well-known that every such cover is normal, meaning that there is a non-trivial deck transformation $\tau:\hat\Sigma\to\hat\Sigma$.  Clearly $\tau$ is an involution, i.e. order $2$. Notice also that associated to any such degree 2 cover there is a non-trivial homomorphism
$$\sigma_p:\pi_1(\Sigma)\to\Z/2\Z$$
determined by $\pi_1(\hat\Sigma)=\Ker(\sigma_p)$.

In fact, we can show as follows that there is a, far from unique, non-separating simple closed curve $\alpha$ with $\sigma_p(\gamma)=\langle\alpha,\gamma\rangle\mod 2$.  Any map $\pi_1(\Sigma)\to\Z/2\Z$ factors through $\HH_1(\Sigma, \Z/2\Z)$.  By \eqref{eq:preli1}, the intersection pairing on $\HH_1(\Sigma, \Z/2\Z)$ is
non-degenerate and any map $$\HH_1(\Sigma, \Z/2\Z) \to \Z/2\Z$$ is equivalent to pairing with some (unique) $\Z/2\Z$ homology class.  By \cite[Prop 6.2]{FarbMarg}, any primitive
$\Z$-homology class and hence any $\Z/2\Z$-homology class is represented by a simple closed curve.  Conversely, every non-separating simple closed curve $\alpha\subset\Sigma$ yields 
a homomorphism $\pi_1(\Sigma)\to\Z/2\Z$ and hence a degree 2 cover.

We will frequently talk of {\em lifts} with respect to the cover $p: \hat{\Sigma} \to \Sigma$.  For a homeomorphism, $\varphi: \Sigma \to \Sigma$, a {\em lift}
$\hat{\varphi}$ of $\varphi$ will mean a homeomorphism $\hat{\varphi}: \hat{\Sigma} \to \hat{\Sigma}$ such that $p \circ \hat{\varphi} = \varphi \circ p$.  Simple closed 
curves on $\Sigma$ will be viewed as embedded, connected $1$-dimensional submanifolds of $\Sigma$.  A {\em lift} of a simple closed curve $\gamma$ is a simple closed curve
$\hat \gamma$ in $\hat \Sigma$ which maps homeomorphically to $\gamma$ via $p$.

\subsection{Mapping classes}
Mapping classes come in three different varieties as mentioned above.  There are {\em pseudo-Anosov} mapping classes which are classes having representatives preserving a pair of transverse projective measured foliations of the surface.    Another kind are {\em reducible} mapping
classes; these preserve an isotopy class of an essential 1-dimensional submanifold, i.e. a disjoint collection of simple closed curves, none of which is homotopically trivial. Finally, there are mapping classes which have finite order and these are called {\em periodic}.  The following theorem \cite[Theorem 13.2]{FarbMarg} says that these are the only possibilities.

\begin{theorem}[Nielsen--Thurston]
All mapping classes are periodic, reducible, or pseudo-Anosov.  Furthermore, a mapping class is pseudo-Anosov if and only if it is neither periodic nor reducible.
\end{theorem}

The Torelli group is torsion-free and hence does not contain elements of finite order other than the identity \cite[Theorem 6.12]{FarbMarg}. It is also known that elements in the Torelli group are {\em pure}, meaning that whenever a mapping class in the Torelli group preserves a $1$-dimensional submanifold, then it actually preserves each component (with orientation) 
\cite[Theorem 1.2]{Ivan}.  Thus, the above theorem amounts to the following for the Torelli group.

\begin{corollary} \label{corollary:ClassTor}
Let $\varphi$ be an element of the Torelli group.  Then $\varphi$ is pseudo-Anosov if and only if it does not fix the isotopy class of an essential oriented simple closed curve.
\end{corollary}

\section{Homomorphisms From Covers} \label{section:defining}
Throughout this section suppose that $\hat\Sigma$ is a connected surface and $p:\hat\Sigma\to\Sigma$ a degree 2 cover. We associate to $p$ a homomorphism from a finite index subgroup $\Gamma_p$ of $\Mod(\hat\Sigma)$ to $\PSp(2g-2,\Z)=\Sp(2g-2,\Z)/\pm\Id$. A similar construction appears for example in \cite{Looijenga}. Our main result 
in this section is that the image of the Torelli group $\Torelli(\Sigma)$ under this homomorphism has finite index in $\PSp(2g-2,\Z)$.

From now on we denote by
$$K_p=\Ker(p_*:\HH_1(\hat\Sigma,\Z)\to\HH_1(\Sigma,\Z))$$
the kernel of the transformation in homology induced by the covering map $p$. Notice that $K_p$ is invariant under the action in homology induced by the non-trivial deck transformation $\tau:\hat\Sigma\to\hat\Sigma$.

\begin{lemma} \label{lemma:kernel}        
There is a genus $g-1$ embedded subsurface $X \subseteq \Sigma$ with one boundary component such that $p^{-1}(X)$ is the disjoint union of two subsurfaces $X_\pm$ of $\hat\Sigma$.
Furthermore, $(\Id - \tau_*): \HH_1(X_+, \Z) \to K_p$ is an isomorphism, and $\tau_*|_{K_p} = -\Id$
\end{lemma}
\begin{proof}
Let $\alpha\subset\Sigma$ be a non-separating curve such that the cover $p:\hat{\Sigma} \to \Sigma$ is determined by a homomorphism 
$$\sigma: \pi_1(\Sigma) \to \Z/2\Z,\ \ \sigma(\eta) = \langle \eta, \alpha \rangle \;\;\; \text{mod } 2$$
and choose a simple closed curve $\beta$ intersecting $\alpha$ exactly once.  The two curves $\alpha$ and $\beta$ fill a one holed torus $T$.  Let $X = \Sigma \setminus \interior{T}$ be the complement of the interior of $T$.

The subgroup $\pi_1(X)$ is contained in the kernel of $\sigma$, and so the preimage of $X$ under $p$ has two components $X_+$ and $X_-$, each of which $p$ identifies with $X$.  The deck transformation $\tau$ interchanges them.  A Mayer-Vietoris argument shows that the following submodule is a direct summand of $\HH_1(\hat\Sigma, \Z)$:
$$\HH_1(X_+ \cup X_-, \Z) \cong \HH_1(X_+, \Z) \oplus \HH_1(X_-, \Z)$$
Since $\tau_*$ exchanges $\HH_1(X_+, \Z)$ and $\HH_1(X_-, \Z)$, the submodule 
$(\Id - \tau_*)\HH_1(X_+, \Z)$ is a direct summand of $\HH_1(X_+, \Z) \oplus \HH_1(X_-, \Z)$ and hence of $\HH_1(\hat\Sigma, \Z)$.  Furthermore, observe 
$(\Id - \tau_*)\HH_1(X_+, \Z)$ is contained in $K_p$.  These facts, combined with considerations of rank, imply that $(\Id - \tau_*)\HH_1(X_+, \Z) = K_p$. 

Noticing that $\tau_*$ is an involution we obtain $\tau_*\circ(\Id - \tau_*)=-(\Id-\tau_*)$. The last claim follows.
\end{proof}

We continue with the same notation. A computation shows that the isomorphism
$$(\Id - \tau_*): \HH_1(X_+, \Z) \to K_p$$
scales the intersection pairing in $\HH_1(X_+, \Z)$ and $K_p$; more precisely we have
$$\langle(\Id - \tau_*)(\cdot),(\Id - \tau_*)(\cdot)\rangle=2\langle\cdot,\cdot\rangle$$
Consequently, the intersection pairing is non-degenerate on $K_p$ and 
$$\Aut(K_p,\langle\cdot,\cdot\rangle) \cong\Aut(\HH_1(X_+,\Z),\langle\cdot,\cdot\rangle)\cong \Sp(2g-2, \Z)$$ 
It follows hence from the last claim of Lemma \ref{lemma:kernel} that
$$\Aut(K_p,\langle\cdot,\cdot\rangle)/\langle\tau_*\rangle \cong \PSp(2g-2,\Z)$$
Our next goal is to associate to the cover $p:\hat\Sigma\to\Sigma$ a homomorphism from a finite index subgroup $\Gamma_p$ of $\Mod(\Sigma)$ to $\PSp(2g-2,\Z)$. 

\begin{lemma} \label{lemma:fphom}
The subgroup $\Gamma_p\subset\Mod(\Sigma)$ consisting of those elements which have a lift to $\hat\Sigma$ has finite index in $\Mod(\Sigma)$ and contains $\Torelli(\Sigma)$.
\end{lemma}
\begin{proof}
As in the proof Lemma \ref{lemma:kernel} consider the homomorphism $\sigma:\pi_1(\Sigma)\to\Z/2\Z$ and recall that $\pi_1(\hat\Sigma)$ is the kernel of $\sigma$. Notice that the normality of $\pi_1(\hat\Sigma)$ implies that there is no ambiguity in considering $\pi_1(\hat\Sigma)$ and its image under a mapping class as subgroups of $\pi_1(\Sigma)$. From this point of view, the group $\Gamma_p$ is precisely the group of mapping classes preserving the subgroup $\pi_1(\hat\Sigma)$. Since $\pi_1(\hat\Sigma)$ has finite index, it follows that $\Gamma_p$ has finite index in $\Mod(\Sigma)$. Given $\phi\in\Torelli(\Sigma)$ we have $\sigma=\sigma\circ\phi$; it follows that $\phi$ preserves $\pi_1(\hat\Sigma)$.  This proves that $\Torelli(\Sigma)\subset\Gamma_p$.
\end{proof}

At this point we can construct the desired homomorphism. Given $\varphi\in\Gamma_p$ let $\hat\varphi\in\Mod(\hat\Sigma)$ be a lift and notice that $\hat\varphi=\tau\circ\hat\varphi\circ\tau$. It follows that $\hat\varphi_*:\HH_1(\hat\Sigma,\Z)\to\HH_1(\hat\Sigma,\Z)$ preserves the submodule $K_p$. In particular, $\hat\varphi_*\vert_{K_p}$ represents an element in $\Aut(K_p,\langle\cdot,\cdot\rangle)$. 

The element $\varphi\in\Gamma_p$ has a second lift, namely $\tau\circ\hat\varphi$. The automorphisms of $K_p$ represented by $(\tau\circ\hat\varphi)_*$ differs from $\hat\varphi_*$ by $\tau_*$. In particular, we obtain a well-defined homomorphism
\begin{equation}\label{lemma:welldefhom}
\ell_p:\Gamma_p\to\Aut(K_p,\langle\cdot,\cdot\rangle)/\langle\tau_*\rangle\cong\PSp(2g-2,\Z)
\end{equation}

Our next goal is to prove the following lemma.

\begin{lemma} \label{lemma:finiteindex}
Let $\Sigma$ be closed surface of genus $g \geq 3$.  The image of $\Torelli(\Sigma)$ in $\PSp(2g-2, \Z)$ under the homomorphism $\ell_p$ has finite index in $\PSp(2g-2, \Z)$.
\end{lemma}

In order to prove Lemma \ref{lemma:finiteindex} we will show that $\Gamma_p$ surjects onto $\PSp(2g-2, \Z)$ and that the image of $\Torelli(\Sigma)$ is infinite; the claim follows then from Margulis's Normal Subgroup Theorem.

\begin{lemma}\label{lem-surject}
The homomorphism $\ell_p: \Gamma_p\to \PSp(2g-2, \Z)$ is surjective.
\end{lemma}

\begin{proof}
Let $X$ be as provided by Lemma \ref{lemma:kernel} and recall that $p^{-1}(X)$ is the disjoint union of two subsurfaces $X_+$ and $X_-$ of $\hat\Sigma$. Notice that $p|_{X_+}$ and $p|_{X_-}$ are both homeomorphisms and let $q: X \to X_+$ be the inverse of $p|_{X_+}$; the map $\tau \circ q$ is the inverse of $p|_{X_-}$. Notice that the homomorphism
\begin{equation}\label{eq-bla}
\iota=(\Id-\tau_*)\circ q_*:\HH_1(X,\Z)\to K_p
\end{equation}
is an isomorphism by Lemma \ref{lemma:kernel}.

Denote by $\Homeo(\Sigma,\Sigma\setminus X)$ the group of those homeomorphisms which fix $\Sigma\setminus X$ pointwise. Given $\varphi\in\Homeo(\Sigma,\Sigma\setminus X)$ we can define $\hat\varphi\in\Homeo(\hat{\Sigma})$ as
$$\begin{array}{lcl} \hat\varphi |_{\hat{\Sigma} \setminus p^{-1}(X)} &=& \Id\\
				\hat\varphi |_{X_+} &=& q \circ\varphi \circ p \\
				\hat\varphi |_{X_-} &=& \tau \circ q \circ\varphi\circ p \end{array}.
$$
By definition, $\hat\varphi$ is a lift of $\varphi$ and thus preserves $K_p$; by construction we have $\hat\varphi_*\vert_{K_p}=\iota\circ\varphi_*\circ\iota^{-1}$ where $\iota$ is as in \eqref{eq-bla}.

Since the homomorphism 
$$\Homeo(\Sigma,\Sigma\setminus X)\to\Aut(H_1(X,\Z),\langle\cdot,\cdot\rangle),\ \ \varphi\mapsto\varphi_*$$
is surjective, it follows that the homomorphism
$$\Homeo(\Sigma,\Sigma\setminus X)\to\Aut(K_p,\langle\cdot,\cdot\rangle),\ \ \varphi\mapsto\hat\varphi_*\vert_{K_p}=\iota\circ\varphi_*\circ\iota^{-1}$$
is also surjective. The induced homomorphism to $\Aut(K_p,\langle\cdot,\cdot\rangle)/\langle\tau_*\rangle\cong\PSp(2g-2,\Z)$ is then a fortiori surjective. We have proved that the restriction of $\ell_p$ to the image of $\Homeo(\Sigma,\Sigma\setminus X)$ in $\Gamma$ is surjective; the claim follows.
\end{proof}

Our next goal is to show that the image under $\ell_p$ of $\Torelli(\Sigma)$ is infinite. In order to do so we will prove that the image of the Dehn twist along some well-chosen separating curve $\gamma\subset\Sigma$ has infinite order. Recall that any such curve determines a symplectic splitting 
$$\HH_1(\Sigma,\Z)=\HH_1(X_1,\Z)\oplus\HH_1(X_2,\Z)$$
where $X_1,X_2$ are the connected components of $\Sigma\setminus\gamma$; this splitting induces also a splitting of  $\HH_1(\Sigma,\Z/2\Z)$. The preimage $p^{-1}(\gamma)$ of every separating curve $\gamma\subset\Sigma$ consists of two components $\hat\gamma_1,\hat\gamma_2\subset\hat\Sigma$ and these two curves are either both separating or both nonseparating.

\begin{lemma} \label{lemma:seplift}
Let $\gamma$ be a separating simple closed curve on $\Sigma$ and $\HH_1(\Sigma,\Z/2\Z)=V_1 \oplus V_2$ the induced splitting.
 The lifts of $\gamma$ to $\hat\Sigma$ are separating if and only if there is some $i=1,2$ with $V_i \subset p_*(\HH_1(\hat{\Sigma}, \Z/2\Z))$.
\end{lemma}
\begin{proof}
Let $X_1$ be the subsurface on one side of $\gamma$ and $X_2$ the one on the other side so that $V_i=\HH_1(X_i, \Z/2\Z)$. Suppose the two lifts $\hat\gamma_1$ and $\hat\gamma_2$ of $\gamma$ to $\hat\Sigma$ are separating. Since $\hat{\Sigma} \setminus (\hat{\gamma}_1\cup\hat{\gamma}_2)$ consists of three connected components $Z_1, Z_2, Z_3$ and $\tau$ is order $2$, it must interchange two components.  Consequently, $p$ induces a homeomorphism from some $Z_j$ to some $X_i$ and this implies one of the $V_i$ lies in
$p_*(\HH_1(\hat{\Sigma}, \Z/2\Z))$.

Conversely, if one of the $V_i$ is contained in $p_*(\HH_1(\hat{\Sigma}, \Z/2\Z))$, then the image of $\pi_1(X_i)$ under the inclusion map lies in the image of 
$\pi_1(\hat{\Sigma})$ and so the inclusion map lifts to the cover.  Since $\gamma = \partial X_i$, it lifts to a bounding curve in the cover.
\end{proof}

\begin{lemma} \label{lemma:nicesepcurve}
For any $2$-fold covering map $p: \hat{\Sigma} \to \Sigma$, there is a separating simple closed curve $\gamma$ on $\Sigma$ with nonseparating lifts.  Conversely, for any
separating simple closed curve $\gamma$ on $\Sigma$, there is a $2$-fold covering map $p: \hat{\Sigma} \to \Sigma$ such that $\gamma$ has nonseparating lifts.
\end{lemma}
\begin{proof}
Since the mod $2$ intersection pairing is nondegenerate, 
codimension $1$ subspaces of $\HH_1(\Sigma, \Z/2\Z) \cong (\Z/2\Z)^{2g}$ correspond bijectively to $1$-dimensional subspaces and hence to non-zero vectors. The surjectivity of $\Mod(\Sigma)\to\Sp(2g,\Z/2\Z)$ implies that $\Mod(\Sigma)$ acts transitively on non-zero vectors of $\HH_1(\Sigma, \Z/2\Z)$ and thus on the set of codimension $1$ subspaces therein.

Now, let $\eta$ be an arbitrary separating simple closed curve and let $p: \hat{\Sigma} \to \Sigma$ be an arbitrary $2$-fold cover.  Let $V_1 \oplus V_2$ be the splitting
of $\HH_1(\Sigma, \Z/2\Z)$ induced by $\eta$, and note that there is some codimension $1$ subspace containing neither $V_1$ nor $V_2$.
The above implies there is some homeomorphism $\varphi$ of $\Sigma$ such that $\varphi_*(p_*(\HH_1(\hat{\Sigma}, \Z/2\Z)))$ contains neither $V_1$ nor $V_2$, and so
by Lemma \ref{lemma:seplift}, $\eta$ has nonseparating lifts with respect to the covering map $\varphi \circ p$.  Conversely, we see that $\varphi^{-1}(\eta)$ has nonseparating
lifts with respect to the covering $p$.
\end{proof}

\begin{lemma} \label{lemma:inforder}
Let $p: \hat{\Sigma} \to \Sigma$ be a degree $2$ cover.  There is a mapping class $\varphi \in \Torelli(\Sigma)$ such that $\ell_p(\varphi)$ has infinite order.
\end{lemma}
\begin{proof}
Let $\gamma$ be as in Lemma \ref{lemma:nicesepcurve} and $\hat\gamma,\tau(\hat\gamma)$ the two components of $p^{-1}(\gamma)$; here $\tau$ is as always the non-trivial deck transformation of the cover. 

Then the right Dehn twist, $T_\gamma$ about $\gamma$ lifts to the composition of the right twists $T_{\hat{\gamma}} T_{\tau(\hat{\gamma})}$. Noting that the classes $\tau(\gamma)$ and $-\gamma$ are equal in $\HH_1(\hat\Sigma,\Z)$ we have for any $v\in\HH_1(\hat\Sigma,\Z)$:
$$T_{\hat{\gamma}} T_{\tau(\hat{\gamma})}(v) = v + 2\langle\hat\gamma,v\rangle\hat{\gamma}$$
Notice that $\gamma\in K_p$ is non-trivial and that there is some $v\in K_p$ with $\langle\gamma,v\rangle\neq 0$ because the intersection pairing is non-degenerate on $K_p$. For any such $v$ we deduce that $T_{\hat{\gamma}}^m T_{\tau(\hat{\gamma})}^m(v)\neq v$ for all $m\neq 0$. This implies that $\ell_p(T_\gamma)$ has infinite order.
\end{proof}

We are now ready to prove Lemma \ref{lemma:finiteindex}:
 
\begin{proof}[Proof of Lemma \ref{lemma:finiteindex}]
Since $\Gamma$ surjects onto $\PSp(2g-2, \Z)$ by Lemma \ref{lem-surject}, $\ell_p(\Torelli(\Sigma))$ is normal.  By Lemma \ref{lemma:inforder}, $\ell_p(\Torelli(\Sigma))\subset\PSp(2g-2,\Z)$ is infinite. Hence, Margulis's Normal Subgroups Theorem \cite{Margulis} implies that $\ell_p(\Torelli(\Sigma))$ has finite index in $\PSp(2g-2,\Z)$.
\end{proof}

The same proof applies mutatis mutandis to the following.
\begin{corollary} \label{corollary:Kgfi}
Let $\Sigma$ be closed surface of genus $g \geq 3$.  The image of ${\mathcal K}(\Sigma)$ in $\PSp(2g-2, \Z)$ under the homomorphism $\ell_p$ has finite index in $\PSp(2g-2, \Z)$.
\end{corollary}

\section{A Criterion for Being Pseudo-Anosov}
In this section we provide a criterion to detect pseudo-Anosov elements in the Torelli group.  This criterion is based on the action of the homology on index $2$ covers.  While it is in 
the spirit of the Casson-Bleiler criterion \cite[Lemma 5.1]{CassBlei}, it is not literally the same. Notice for instance that the assumptions in the original Casson-Bleiler criterion are not satisfied by any lift of a mapping class.  Indeed, as we have seen, any lift necessarily preserves a proper subspace of homology and so will never have irreducible characteristic polynomial.

The criterion works roughly as follows.  By Corollary \ref{corollary:ClassTor}, a reducible mapping class $\varphi$ in the Torelli group fixes a simple closed curve.  We will show that there is a 2-sheeted cover $p:\hat\Sigma\to\Sigma$ such that $\ell_p(\varphi)\in\PSp(2g-2,\Z)$ fixes a line in $\Z^{2g-2}$. Equivalently, any lift of $\ell_p(\varphi)$ to $\Sp(2g-2,\Z)$ has an eigenvector in $\Q^{2g-2}$. We prove:

\begin{theorem} \label{theorem:criterion}
Let $\varphi\in \Torelli(\Sigma)$ and suppose that for every 2-sheeted cover $p:\hat\Sigma\to\Sigma$ the element $\ell_p(\varphi)$ has no invariant line in $\Z^{2g-2}$. Then $\varphi$ is  pseudo-Anosov.  
\end{theorem}

\begin{proof}
Suppose that $\varphi$ is not pseudo-Anosov and recall that by Corollary \ref{corollary:ClassTor} the mapping class $\varphi$ fixes a simple (oriented) closed curve $\gamma$ in $\Sigma$. We will show that there is a 2-sheeted cover $p:\hat\Sigma\to\Sigma$ such that $\ell_p(\varphi)$ fixes a line in $\Z^{2g-2}$.  The proof breaks into two cases: $\gamma$
is separating and $\gamma$ is nonseparating.  
\medskip

\paragraph{\bf The curve $\gamma$ is separating.} By Lemma \ref {lemma:nicesepcurve}, there is a 2-sheeted cover $p:\hat\Sigma\to\Sigma$ such that the lifts $\hat\gamma,\tau(\hat\gamma)$ of $\gamma$ to $\hat\Sigma$ are nonseparating. Here $\tau$ is the non-trivial deck transformation of $p$. Notice that the class $\hat\gamma\in\HH_1(\hat\Sigma,\Z)$ is non-trivial, that $\hat\gamma\in K_p$ and $\tau(\hat\gamma)=-\hat\gamma$. Any lift $\hat\varphi$ of $\varphi$ to $\hat\Sigma$ either fixes the curves $\hat\gamma,\tau(\hat\gamma)$ or permutes them. It follows that $\hat\varphi$ fixes the line in $K_p$ spanned by $\hat\gamma$. We are done in this case.
\medskip

\paragraph{\bf The curve $\gamma$ is nonseparating.} Let $\alpha\subset\Sigma$ be a non-separating curve representing a homology class with $\langle\alpha,\gamma\rangle=0$ but such that the classes $\alpha$ and $\gamma$ are distinct in $\HH_1(\Sigma,\Z/2\Z)$. Consider the homomorphism
$$\pi_1(\Sigma)\to\Z/2\Z,\ \ \eta\mapsto\langle\eta,\alpha\rangle \;\; \text{mod } 2$$
and the associated 2-sheeted cover
$$p:\hat\Sigma\to\Sigma$$
Since $\langle\alpha,\gamma\rangle=0$, we have that $p^{-1}(\gamma)$ consists of two components $\hat\gamma,\tau(\hat\gamma)$ where $\tau$ is as always the non-trivial deck transformation of the cover in question. We claim that $\hat\gamma-\tau(\hat\gamma)$ is a non-trivial homology class. Suppose that this is not the case. Then $p^{-1}(\gamma)$ separates $\hat{\Sigma}$.  Consequently, any closed curve (transversal to $p^{-1}(\gamma)$) must intersect $p^{-1}(\gamma)$ an even number of times.  Thus, all elements in $p_*(\pi_1(\hat{\Sigma}))$ have $0$ mod $2$ intersection with $\gamma$.  Since $p_*(\pi_1(\hat{\Sigma}))$ is index $2$ in $\pi_1(\Sigma)$, it must consist of all elements which have $0$ mod $2$ intersection with 
$\gamma$. It follows that $\alpha$ and $\gamma$ are equal in $\HH_1(\Sigma,\Z/2\Z)$ contradicting our assumption.

At this point we know that $\hat\gamma-\tau(\hat\gamma)$ is non-trivial in $\HH_1(\hat\Sigma,\Z)$. Observing again that any lift $\hat\varphi$ of $\varphi$ to $\hat\Sigma$ either fixes the curves $\hat\gamma$ and $\tau(\hat\gamma)$ or permutes them, it follows that $\ell_p(\varphi)$ fixes the line in $K_p$ spanned by $\hat\gamma-\tau(\hat\gamma)$.
\end{proof}

Theorem \ref{theorem:criterion} implies directly the following useful fact:

\begin{corollary}\label{kor-useful}
Let $\varphi\in \Torelli(\Sigma)$ and suppose that for every 2-sheeted cover $p:\hat\Sigma\to\Sigma$ there is a matrix in $\Sp(2g-2,\Z)$ representing $\ell_p(\varphi)\in\PSp(2g-2,\Z)$ which has irreducible characteristic polynomial. Then $\varphi$ is pseudo-Anosov.
\end{corollary}

\section{Genericity of Pseudo-Anosovs}
We now combine the results from the previous sections with theorems from \cite{Rivin} to prove our main theorem.  We restate one of the results from \cite{Rivin} in a slightly more convenient way.
\begin{theorem}[Rivin] \label{theorem:Rivin}
Let $S\subset\Sp(2g,\Z)$ be a symmetric finite generating set of a Zariski dense subgroup $H < \Sp(2g, \C)$.  Let $w \in W_n(S)$ be a word chosen uniformly at random and $P_n$ be the probability that $h(w)$
has irreducible characteristic polynomial.  Then, $1-P_n = O(c^n)$ for some $c < 1$.
\end{theorem}
\begin{proof}
Combine \cite[Theorem 3.1]{Rivin} and \cite[Theorem 2.4]{Rivin}.
\end{proof}
Before proving our main theorem, notice that any finite index subgroup of $\Sp(2g,\Z)$ is a lattice in $\Sp(2g,\R)$ and hence Zariski dense in $\Sp(2g,\C)$ by Borel's Density Theorem; see for example \cite[Section 4.7]{Witt}. We are now ready to prove Theorem \ref{theorem:maintheorem}:

\begin{named}{Theorem \ref{theorem:maintheorem}}
Let $\Sigma$ be a closed surface of genus $g \geq 3$, let $S$ a finite symmetric generating set of $\Torelli(\Sigma)$, and denote by $P_n$ the probability that the element $h(w)\in\Torelli(\Sigma)$ represented by a word $w \in W_n(S)$ chosen uniformly at random is pseudo-Anosov.  Then, $1- P_n = O(c^n)$ for some $c < 1$.
\end{named}
\begin{proof}
A surface $\Sigma$ of genus $g$ has $2^{2g}-1$ covers of degree 2. Consider for each such cover the associated homomorphism 
$$\ell_p:\Torelli(\Sigma)\to\PSp(2g-2,\Z)$$
Let us say (for the purposes of this proof) that $\ell_p(\varphi)$ is {\em irreducible} if a representative matrix of $\ell_p(\varphi)$ in $\Sp(2g-2, \Z)$ has irreducible characteristic polynomial and is {\em  reducible}
otherwise. 
Observe that, since there are only finitely many covers $p$, it suffices to show by Corollary \ref{kor-useful} that the probability that $h(w)$ for $w \in W_n(\ell_p(S))$ is reducible is $O(c^n)$ for some 
$c < 1$.

We now need to translate our situation to $\Sp$ from $\PSp$ to apply Theorem \ref{theorem:Rivin}.  Recall that $W_n(S)$ is the set of words of length $n$ in a set of generators $S$, and $h(w)$ is the group element corresponding to $w \in W_n(S)$.  Also let $$\pi: \Sp(2g-2, \Z) \to \PSp(2g-2, \Z)$$ be the quotient map, and  
let $S_p = \pi^{-1}(\ell_p(S))$.  Notice that for every $\ell_p(s) \in \ell_p(S)$, the set $\pi^{-1}(\{\ell_p(s)\})$ has two elements, say $s_p, s'_p$, and $s_p^{-1} s'_p = -\Id$ 
which generates the kernel of $\pi$.  Consequently, $S_p$ generates $\pi^{-1}(\ell_p(\Torelli(\Sigma)))$ which is a finite index subgroup of $\Sp(2g-2, \Z)$.

Suppose now that $w$ is an arbitrary word of length $n$ in the letters $\ell_p(S)$.  Consider the natural projections of words of length $n$ (not of group elements) 
$W_n(S_p) \to W_n(\ell_p(S))$.  Since every $\ell_p(s) \in \ell_p(S)$ has exactly two preimages
in $S_p$, there are exactly $2^n$ words in $W_n(S_p)$ that project to $w$.  Consequently, for any group element $g \in \PSp(2g-2, \Z)$, the following probabilities are equal.
\begin{equation} P(h(w) = g \; | \; w \in W_n(\ell_p(S))) = P(\pi(h(w)) = g \;  | \; w \in W_n(S_p)) \label{eqn:prob} \end{equation}
Now by Theorem \ref{theorem:Rivin} and \eqref{eqn:prob}, 
$$P(h(w) \text{ is reducible} \; |\; w \in W_n(\ell_p(S))) = O(c^n).$$
\end{proof}

\section{Further Results}
In this section we prove Theorem \ref{theorem:maintheorem2}.  Its proof is the same as that for Theorem \ref{theorem:maintheorem} with some minor additions which we present now.  Since general mapping classes are not pure, we must prove a new version of Theorem \ref{theorem:criterion}.

\begin{corollary}[Corollary to proof of Theorem \ref{theorem:criterion}]
Let $\varphi \in  \cap_p \Gamma_p$.  If for every 2-sheeted cover $p: \hat \Sigma \to \Sigma$, some representative matrix in $\Sp(2g-2, \Z)$ of the element $\ell_p(\varphi)$ has
no roots of unity as an eigenvalue, then $\varphi$ is pseudo-Anosov.
\end{corollary}

\begin{proof}
Suppose $\varphi$ is not pseudo-Anosov.  Then it is periodic or reducible, and this implies that some power $\varphi^n$ fixes a simple closed curve.  By the proof of Theorem \ref{theorem:criterion}, for some $p$, 
the map $\ell_p(\varphi^n)$ fixes a line in $\Z^{2g-2}$, and so for any representative matrix $M < \Sp(2g-2, \Z)$ of $\ell_p(\varphi)$, it must be that $M^n$ has eigenvalue $\pm 1$.
Hence, $M^{2n}$ has eigenvalue $1$, and this implies $M$ has some root of unity as an eigenvalue.
\end{proof} 

We now prove the following.
\begin{named}{Theorem \ref{theorem:maintheorem2}}
Let $\Sigma$ be a surface of genus $g > 2$.  Let $G$ be a finitely generated subgroup of $\cap_p \Gamma_p$ with symmetric generating set $S$, and let $G_p < \Sp(2g-2, \Z)$ be the full preimage of $\ell_p(G) < \PSp(2g-2, \Z)$.  If $G_p$ is Zariski dense for all index $2$ covers $p$, 
then the probability of a random walk in $S$ being pseudo-Anosov is $1-O(c^n)$ for some $c < 1$.
\end{named}
\begin{proof}
Let $S_p \subset \Sp(2g-2, \Z)$ be the full preimage of $\ell_p(S) \subset \PSp(2g-2, \Z)$.  As in the proof of Theorem 
\ref{theorem:maintheorem}, we are reduced to showing that, away from an exponentially small set, a word in $S_p$ yields a matrix that has no roots of unity as eigenvalues.  Observe that $S_p$ generates the full preimage of $\ell_p(G)$ which is 
Zariski dense by assumption.

Let $w$ be a random word of length $n$ in $S_p$ representing an element $h(w) \in \Sp(2g-2, \Z)$ with characteristic polynomial $q(x)$.  By \cite[Theorem 7.12]{Kowalski}, with probability $1 - O(c^n)$, the polynomial $q(x)$ is irreducible and has Galois group isomorphic to the signed permutation group on $g-1$ letters which has order $2^{g-1} (g-1)!$.  This will imply that 
$h(w)$ has no roots of unity as an eigenvalue.

Suppose $h(w)$ did have a root of unity $\zeta$ as an eigenvalue.  Then, $\zeta$ is a primitive $n$th root of unity for some $n$.  Then since $q(x)$ has $\Q$-coefficients, the 
minimal polynomial for $\zeta$ over $\Q$, namely the $n$th cyclotomic polynomial, must divide $q(x)$.  Since $q(x)$ is irreducible, it must be the cyclotomic polynomial.  However, the Galois group of a cyclotomic polynomial
has order equal to its degree which is $\text{deg } q(x) = 2g-2$. (See Chapter VI \S 3 of \cite{Lang} for these facts about cyclotomic polynomials.)  Since $g > 2$, the order of the Galois group gives a contradiction.  

\end{proof}

Corollary \ref{cor:maincor} now follows from Corollary \ref{corollary:Kgfi}, Theorem \ref{theorem:maintheorem2}, and the fact that finite index subgroups of $\Sp(2g-2, \Z)$ are 
Zariski dense.


\begin{thebibliography}{0}

\bibitem{CassBlei} Casson, A., and S. Bleiler. 
\textit{Automorphisms of surfaces after Nielsen and Thurston},
Cambridge: Cambridge University Press, 1988.

\bibitem{FarbMarg} Farb, B., and D. Margalit.
\textit{A primer on mapping class groups},
Princeton, NJ: Princeton University Press, 2011.

\bibitem{Ivan} Ivanov, N. 
\textit{Subgroups of Teichm\"uller modular groups},
Translations of Mathematical Monographs, vol. 115,
Providence, RI: American Mathematical Society, 1992.

\bibitem{Johnson}
Johnson, D.
``The structure of the Torelli group. I. A finite set of generators for $\Torelli$.''
\textit{Annals of Mathematics (2)} 118, no. 3 (1983): 423--442.

\bibitem{Kowalski}
Kowalski, E.
\textit{The large sieve and its applications. Arithmetic geometry, random walks and discrete group},
Cambridge Tracts in Mathematics 175. 
Cambridge: Cambridge University Press, 2008.

\bibitem{Lang}
Lang, S.
{\em Algebra}, 
Graduate Texts in Mathematics 211.
New York: Springer-Verlag, 2002.

\bibitem{Looijenga}
Looijenga, E.
``Prym representations of mapping class groups.''
\textit{Geometriae Dedicata} 64, no. 1 (1997): 69Ð83. 

\bibitem{Lubotzky}
Lubotzky A., and C. Meiri.
``Sieve methods in group theory: II, The Torelli subgroup.'' preprint.

\bibitem{Maher} 
Maher, J. 
``Random walks on the mapping class group.'' 
\textit{Duke Mathematical Journal} 156, no. 3 (2011):  429--468.

\bibitem{Margulis}
Margulis, G., 
``Finiteness of quotient groups of discrete subgroups.''
\textit{Functional Analysis and its Applications} 13 (1979): 178-187.

\bibitem{MM} 
McCullough, D., and A. Miller.
``The genus $2$ Torelli group is not finitely generated.''
{\em Topology and its Applications} 22, no. 1 (1986): 43--49.

\bibitem{Witt} 
Morris, D. W.
\textit{Ratner's theorems on unipotent flows},
Chicago Lectures in Mathematics.
Chicago, IL: University of Chicago Press, 2005.

\bibitem{Riv2} 
Rivin, I.
``Walks on graphs and lattices--effective bounds and applications.''
\textit{Forum Mathematicum} 21, no. 4 (2009): 673--685.

\bibitem{Riv} 
Rivin, I.
``Walks on groups, counting reducible matrices, polynomials, and surface and free group automorphisms.''
\textit{Duke Mathematical Journal} 142, no. 2 (2008): 353--379. 

\bibitem{Rivin} 
Rivin, I.
``Zariski density and genericity.''
\textit{International Mathematics Research Notices}, no. 19 (2010): 3649--3657. 

\end{thebibliography}
\end{document}